\documentclass[a4paper,10pt]{amsart}
\usepackage[latin1]{inputenc}
\usepackage{amsmath,amsfonts,amssymb}
\def\R{\mathbb{R}}
\def\C{\mathbb{C}}
\def\B{\mathbb{B}}

\def\PP{\mathbb{P}}
\def\ccc{\mathbb{C}}
\def\tsp{{\,}^t\!}
\def\vp{\varphi}
\def\D{\Delta}
\def\bd{\partial\Delta}
\def\pri{\,'\!}

\theoremstyle{plain}
\newtheorem{theorem}{Theorem}
\newtheorem{cor}{Corollary}
\newtheorem{prop}{Proposition}
\newtheorem{lemma}{Lemma}
\newtheorem{definition}{Definition}

\theoremstyle{remark}
\newtheorem{remark}{Remark}

\title{Stationary discs glued to a Levi non-degenerate hypersurface}
\author{L\'ea Blanc-Centi}
\address{L.A.T.P.\\ C.M.I.\\ 39 rue Joliot-Curie\\ 13453 Marseille Cedex 13\\ FRANCE}
\email{lea@cmi.univ-mrs.fr}

\begin{document}

\begin{abstract}
We obtain an explicit parametrization of stationary discs glued to some Levi non-degenerate hypersurfaces. These discs form a family which is invariant under the action of biholomorphisms. We use this parametrization to construct a local circular representation of these hypersurfaces. As a corollary, we get the uniqueness of biholomorphisms with given 1-jet at some convenient point.
\end{abstract}

\maketitle

A central problem in complex analysis consists in classifying the domains in $\C^n$ under the action of biholomorphisms. If $n=1$, the Riemann mapping theorem states that every simply connected domain of $\C$, which is not all of $\C$, is biholomorphic to the unit disc. We know since the work of H. Poincar\'e \cite{Poincare} that this theorem has no immediate generalization in the case $n\ge 2$. Indeed, the rigidity of the holomorphy condition in the multidimensional situation makes biholomorphic equivalence very rare. This leads us to seek biholomorphic invariants of the boundary of a domain. S.S. Chern and J.K. Moser \cite{CM} have associated to any Levi non-degenerate real hypersurface a simple equation and proposed a purely geometric construction of classifying invariants.
Stationary discs are another natural invariant of manifolds with boundary with respect to biholomorphisms, or more generally with respect to CR maps. 

We recall that a holomorphic disc in $\C^n$ is a holomorphic map from the unit disc $\Delta$ in $\C$ to $M$, continuous up to the boundary. L. Lempert \cite{L1} showed that some of these discs have interesting extremality properties with respect to the Kobayashi metric. Indeed, in a strongly convex domain in $\C^n$, the geodesics for this metric are exactly the {\em stationary discs}.
Using this particular family of discs, L. Lempert constructed a homeomorphism between any strongly convex domain and the unit ball, which is an analogue of the Riemann map. More generally, the parametrization of $M$ by the stationary discs yields to a (local) representation of $M$ as a circular submanifold. 
In higher codimension, A. Sukhov et A. Tumanov \cite{ST} obtained, using the same method, a new invariant of small perturbations of the product of spheres $S^3\times S^3$, canonically diffeomorphic to the conormal bundle. The map constructed commutes with biholomorphisms, and is a partial analogue of the circular representation given by L. Lempert.\\

The aim of this paper is to give some explicit parametrization of the stationary discs glued to small deformations of the hyperquadric $Q$ defined by
$$0=r(z)=\Re z_0-\tsp\pri\bar{z}A\pri z,$$ 
where $A$ denotes a $n$-sized non-degenerate hermitian matrix. More precisely: 

\begin{theorem}\label{2diffeo}

Pick $p=(p_0,0,\hdots,0)\notin Q$ and let $h^0$ be a stationary disc glued to $Q$ and centered at $p$.
Then, for any $\rho$ near $r$ with respect to the $\mathcal{C}^3$-topology, the set of all non-constant stationary discs glued to the hypersurface $\{\rho=0\}$ forms a $(4n+3)$-parameter family near $h^0$. 

If we only consider the discs $h$ centered at $p$, the maps $h\mapsto h(1)$ and $h\mapsto h'(0)$ are local diffeomorphisms on their images.
\end{theorem}

Note that Theorem \ref{2diffeo} shows that a hypersurface ``near'' $Q$ is represented in a neighborhood of $p$ in a circular way by its Kobayashi indicatrix $\{h'(0)\}$. The idea of the proof (see \cite{C,ST,CGS-Riemann}) consists in using a criterion of J. Globevnik \cite{Glob}: given a submanifold $E$ and a disc $f$ glued to $E$, and under the condition that the {\em partial indices} of $E$ along $f$ are non-negative, the holomorphic discs near $f$ glued to a submanifold near $E$ form a $\kappa$-parameter family, where $\kappa$ is the {\em Maslov index}. In our case, the main difficulty is to check that partial indices are non-negative, and to compute $\kappa$.\\

Using the biholomorphic invariance property of stationary discs, we get as a corollary that a biholomorphism fixing $p$ in a neighbourhood of small deformations of $Q$ is (locally) uniquely determined by its differential at $p$.\\

This paper is organized as follows. We begin with some definitions in the first section. The second section is devoted to the proof of Theorem \ref{2diffeo}. We then apply this result to the study of biholomorphisms in the third section (Theorem \ref{theo:unicite}). In the Appendices, we prove two quite technical lemmas, used in Section 2. 
\bigskip

\section{Preliminaries}\label{s:Prelim}
Let $M$ be a real hypersurface in $\C^n$ and $\rho$ be some defining function of $M$ (that is, $M=\{\rho=0\}$ and $d\rho$ does not vanish on $M$). A {\em holomorphic disc} $h$ is a holomorphic function from the unit disc $\Delta$ in $\C$ to $\C^n$, continuous up to the boundary. We say that $h$ is {\em glued to $M$} if $h(\partial\Delta)\subset M$.

We recall the construction of the conormal bundle $N^*M$ of $M$. Its fibre at some point $p\in M$ is the set of all $(1,0)$-forms on $T\C^n$ whose real part vanishes identically on the tangent space $T_pM$:
$$N_p^*M=\{\phi\in T_p^*\C^n/\ \Re\phi_{|T_pM}=0\}.$$
Notice that $N_p^*M$ is a (real) line generated by $\partial_p\rho$. We will need the

\begin{prop}\label{prop:conormal}
{\bf \cite{Tumanov}}
A real hypersurface $M$ in $\C^n$ is Levi non-degenerate if and only if its conormal bundle $N^*M$ is totally real out of the zero section.
\end{prop}

\begin{definition}\label{def:stationary}
A holomorphic disc $h$ glued to $M$ is {\em stationary} if there exists a meromorphic lift $(h,h^*)$ of $h$ to the cotangent bundle $T^*\C^n$, with at most one pole of order 1 at 0, such that 
$$\forall\zeta\in\partial\D,\ h^*(\zeta)\in N^*_{h(\zeta)}M\setminus\{0\}.$$
\end{definition}

\begin{remark}\label{rem:defc}
In coordinates, Definition \ref{def:stationary} is equivalent to the existence of some continuous function $c:\partial\D\to\R^*$ such that $h^*(\zeta)=c(\zeta)\partial\rho_{h(\zeta)}$ on $\partial\D$ and $\zeta h^*$ extends holomorphically to $\D$.
This notion is preserved under the action of biholomorphisms.
\end{remark}

Definition \ref{def:stationary} is in fact slightly different from the one given by L. Lempert, since we assume that the boundary of the lift does not meet the zero section of the conormal bundle. In \cite{L1}, even if this condition is not assumed in the definition, it is necessarily satisfied because of the strong convexity of the domain.

Notice that, in light of Proposition \ref{prop:conormal}, the regularity results known for holomorphic discs glued to a totally real submanifold \cite{Chirka} apply to the map $(h,h^*)$. Hence the lift is as regular to the boundary as $N^*M$.
For this reason, following \cite{ST}, we will call in the sequel $(h,h^*)$ (or $h^*$, with a minor abuse of notation) a {\em regular lift of $h$}.\\

A very simple  computation shows that stationary discs glued to the unit sphere and centered at 0 are exactly the linear discs. This result is also valid for general circular strongly convex domains. Indeed, for circular domains, the linear discs are stationary, and conversely for circular strongly convex domains, all stationary discs through 0 are linear by Lempert's uniqueness result. This provides two parametrizations of such discs, \textsl{via} the maps $h\mapsto h(1)$ and $h\mapsto h'(0)$. In the next section, we obtain the same parametrizations only supposing the non-degeneracy of $M$.

\section{Proof of Theorem \ref{2diffeo}}\label{s:Proof}
Consider the hyperquadric $Q\subset\C^{n+1}$ defined on an open set $\Omega\ni 0$ by  
$$0=r(z)=\mathrm{Re}z_0-\sum_{i,j= 1}^na_{i,j}\bar{z}_iz_j,$$
where the hermitian matrix $A=(a_{i,j})_{i,j}$ is non-degenerate. 

\subsection{Stationary discs glued to $Q$}\label{subs:StationaryDiscsQ}
We first determine explicitely stationary discs glued to the hyperquadric.
\begin{prop}
Stationary discs glued to $Q$ are exactly under the form
$$h(\zeta)=\left(\tsp\bar{v}Av+2\tsp\bar{v}Aw\,\frac{\zeta}{1-a\zeta}+\frac{\tsp\bar{w}Aw}{1-|a|^2}\,\frac{1+a\zeta}{1-a\zeta}+iy_0,v+w\,\frac{\zeta}{1-a\zeta}\right)$$
where $v,w\in\C^n$, $y_0\in\R$, $a\in\Delta$.\\
Moreover $h^*$ is a regular lift of $h$ if and only if there exists $b\in\R^*$ such that
$$\forall\zeta\in\Delta\setminus\{0\},\ h^*(\zeta)=\frac{b}{\zeta}\left(\frac{-\bar{a}}{1+|a|^2}+\zeta-\frac{a}{1+|a|^2}\zeta^2\right)\times(1/2,-\tsp\overline{h_\alpha(\zeta)}A).$$
\end{prop}

\begin{remark}\label{exprcentre}
If the disc $h$ is centered at some point $p=(p_0,0,\hdots,0)$, we obtain
$$h(\zeta)=\left(\frac{\tsp\bar{w}Aw}{1-|a|^2}\,\frac{1+a\zeta}{1-a\zeta}+iy_0,w\frac{\zeta}{1-a\zeta}\right).$$  
\end{remark}

\begin{proof}
\underline{Necessary condition}\\
Suppose that $h$ has some regular lift $h^*$, and pick $c$ as in Remark \ref{rem:defc}:
\begin{equation}\label{lienh*h}
\forall\zeta\in\partial\Delta,\ h^*(\zeta)=c(\zeta)\times\frac{\partial r}{\partial z}\circ h(\zeta)=c(\zeta)\times(1/2,-\tsp\overline{h_\alpha(\zeta)}\cdot A).
\end{equation}
Thus, there exists some holomorphic function $\vp$ on $\D$, continuous up to $\partial\D$ and such that for all $\zeta\in\bd$, $c(\zeta)=\vp(\zeta)/\zeta\in\R^*$. Expanding $\vp_{|\partial\D}$ into Fourier series, we obtain that $\vp(\zeta)=a+b\zeta+\bar{a}\zeta^2$ for some $a,b\in\C$, and
\begin{equation}\label{exprc}
\forall\zeta\in\bd,\ c(\zeta)=a\bar{\zeta}+b+\bar{a}\zeta.
\end{equation}
\textbullet\ First case: $a=0$ (and hence $b\not=0$).\\
That is, for all $\zeta\in\bd$, $h^*(\zeta)=b(1/2,-\tsp\overline{h_\alpha(\zeta)}\cdot A)$. By hypothesis, $h^*$ has at most one pole of order 1 at 0, hence $\zeta\bar{h}_\alpha$ is holomorphic in $\D$. Then there exist $v,w\in\C^n$ such that $h_\alpha(\zeta)=v+\zeta w$. Since $h$ is glued to $Q$, we get by means of an expansion into Fourier series that
 $h_0$ is affine, and uniquely defined up to the addition of some purely imaginary constant.\vspace*{0.15cm}\\
\textbullet\ Second case: $a\not=0$.\\
Let us denote by $a_1$ and $a_2$ the zeroes of $0=a+b\zeta+\bar{a}\zeta^2$. Since $h^*$ does not vanish on $\bd$, the moduli of $a_1$ and $a_2$ are different from 1, and $|a_1a_2|=|a/\bar{a}|=1$: assume for example that $0<|a_1|<1<|a_2|$. 
Writing $h_\alpha(\zeta)=\sum_{k=0}^{+\infty} H_k\zeta^k$ in $L^2(\bd)$, we get that the map $\zeta h^*$ extends holomorphically in $\D$ if and only if  
$\zeta\mapsto(a+b\zeta+\bar{a}\zeta^2)\,\overline{h_\alpha(\zeta)}$
extends holomorphically in $\D$, that is,  
\begin{equation}
\forall k\ge 1,\ a\bar{H}_k+b\bar{H}_{k+1}+\bar{a}\bar{H}_{k+2}=0.
\end{equation}
Hence there exist $V, W\in\C^n$, depending only on $H_1$ and $H_2$, such that for all $k\ge 1$, $H_k=\bar{a}_1^{k-1}V+\bar{a}_2^{k-1}W$.
In fact, it is easy to prove that $W=0$, by using the fact that the radius of convergence of each series $\sum (v_j\bar{a}_1^{k-1}+w_j\bar{a}_2^{k-1})\zeta^k$ is larger than 1.
Hence $H_0=h_\alpha(0)$, $H_1=h'_\alpha(0)$ and
\begin{equation}\label{exprhalpha}
h_\alpha(\zeta)=H_0+\zeta H_1\sum_{k=0}^{+\infty}(\bar{a}_1\zeta)^k.
\end{equation}
Since $a_1$ and $a_2$ are not of modulus 1, we get $b\not=0$. By multiplying $h^*$ by $1/b$, we can assume $b=1$. The other regular lifts will be obtained by multiplication by some non-zero constant. We search to express $a$ in terms of $a_1$.
If $1-4|a|^2\le 0$, one gets $|a_1|=|a_2|$, which is impossible.
Consequently, $1-4|a|^2>0$ and
$\displaystyle a_1=\frac{-1+\sqrt{1-4|a|^2}}{2\bar{a}}$. Setting $a=|a|e^{i\theta}$, we obtain 
 $\mathrm{Arg}(a)=\mathrm{Arg}(a_1)-\pi$ and
$|a|=\frac{|a_1|}{1+|a_1|^2}$. Finally,
$\displaystyle a=\frac{-a_1}{1+|a_1|^2}$.\\

\noindent\underline{Sufficient condition}\\ 
Conversely, suppose that $h$ is given as in (\ref{exprhalpha}) 
 (with $|a_1|<1$)
and $h_0$ is uniquely determined (up to the addition of some imaginary constant) in order to glue $h$ to $Q$.
If $a_1=0$, $h_\alpha$ is affine, and the expressions (\ref{lienh*h}) and (\ref{exprc}) show that $h$ is stationary and give its regular lifts.
So let us assume that $0<|a_1|<1$, and define $c(\zeta)=-\frac{a_1}{1+|a_1|^2}\bar{\zeta}+1-\frac{\bar{a}_1}{1+|a_1|^2}\zeta$. Then $h^*$ given by (\ref{lienh*h}) is the only regular lift of $h$ up to multiplicative non-zero constant.
\end{proof}

\begin{prop}\label{forme}
The map $\Phi:(y_0,v,w,a)\mapsto h$ is a smooth diffeomorphism from $\R\times\C^n\times(\C^n\setminus\{0\})\times\D$ to the set of all non-constant stationary discs glued to $Q$. Its inverse is given by 
\begin{eqnarray}\label{parametrisationexplicite}
\Phi^{-1}:h\mapsto \left(\Im h_0(0),h_\alpha(0),h'_\alpha(0), [\theta_h(1)-i\theta_h(i)]\right)
\end{eqnarray}
where
$\displaystyle\theta_h(\zeta)=\frac{||h'_\alpha(0)||^2}{4}\,\left(\frac{1}{||h_\alpha(-\zeta)-h_\alpha(0)||^2}-\frac{1}{||h_\alpha(\zeta)-h_\alpha(0)||^2}\right)$.
\end{prop}

\begin{proof}
By construction, the map $\Phi$ is onto and smooth with respect to the topology induced by $\mathcal{C}^\alpha(\bd)^{n+1}$. Moreover, if $h=\Phi(y_0,v,w,a)$, then $y_0=\Im h_0(0)$ and for all $\zeta\in\D$, $h_\alpha(\zeta)=v+w\,\left[\zeta\sum_{k=0}^{+\infty}(a\zeta)^k\right]$. The uniqueness of the development in power series shows that $\Phi$ is injective.

Let us notice that if $h=\Phi(y_0,v,w,a)$, then for all $\zeta\in\bd$,
$$\frac{||h'_\alpha(0)||^2}{||h_\alpha(\zeta)-h_\alpha(0)||^2}=\left|\frac{1-a\zeta}{\zeta}\right|^2=1-2\Re(\zeta a)+|a|^2,$$
which gives $\Phi^{-1}$. In order to get the smoothness of $\Phi^{-1}$, it suffices to verify that the linear maps $h\mapsto h_\alpha(0)$ and $h\mapsto h'_\alpha(0)$ are continuous on $\mathcal{C}^\alpha(\bar{\D})\cap H(\D)$. But this comes from Cauchy inequalities.

It remains to prove that $\Phi$ is a submersion. As the image of $d\Phi$ is of finite dimension at any point, it suffices to prove that $d\Phi_{(y_0,v,w,a)}$ is injective for all $(y_0,v,w,a)\in\R\times\C^n\times(\C^n\setminus\{0\})\times\D$. Let $(y_0',v',w',a')\in\R\times\C^n\times\C^n\times\C$ be such that $d\Phi_{(y_0,v,w,a)}(y_0',v',w',a')=0$. Writing the partial differentials, we get: 
$$0=d_{y_0}\Phi_{(y_0,v,w,a)}y'_0+d_{v}\Phi_{(y_0,v,w,a)}v'+d_{w}\Phi_{(y_0,v,w,a)}w'+d_{a}\Phi_{(y_0,v,w,a)}a'.$$
The $n$ last components of this equality give
$$\forall \zeta\in\bar{\D},\ v'+w'\frac{\zeta}{1-a\zeta}+w\frac{\zeta^2}{(1-a\zeta)^2}a'=0.$$
Hence $v'=w'=0$ and $wa'=0$, which implies $a'=0$. If we replace in the first component, we obtain  $y'_0=0$, which concludes the proof.
\end{proof}

\begin{remark}
The set of all non-constant stationary discs glued to $Q$ and centered at $p$ is either empty, or a submanifold of (real) dimension $2n+1$.
\end{remark}

\begin{cor}\label{CNSexistence}
Pick $p\in\C^{n+1}$. The set $\mathcal{M}_p$ of all non-constant stationary discs glued to the hyperquadric $Q=\{r=0\}$ and centered at $p$ is non-empty if and only if one of the following conditions holds:

* $A$ is positive definite and $r(p)<0$;

* $A$ is negative definite and $r(p)>0$;

* $A$ has two eigenvalues of different signs.\\
If so, the map $\mathcal{M}_p\ni h\mapsto h(1)\in Q$ is a local diffeomorphism if and only if $p\notin Q$.
\end{cor}

\begin{proof}
Let $h$ be a stationary disc glued to $Q$:
$$h(\zeta)=\left(\tsp\bar{v}Av+2\tsp\bar{v}Aw\,\frac{\zeta}{1-a\zeta}+\frac{\tsp\bar{w}Aw}{1-|a|^2}\,\frac{1+a\zeta}{1-a\zeta}+iy_0,v+w\,\frac{\zeta}{1-a\zeta}\right)$$
where $y_0\in\R$, $v,w\in\C^n$, $a\in\Delta$. In particular, 
\begin{eqnarray}\label{CNScentre}
h(0)=\left(\tsp\bar{v}Av+\frac{\tsp\bar{w}Aw}{1-|a|^2}+iy_0,v\right)\in Q\Leftrightarrow \tsp\bar{w}Aw=0.
\end{eqnarray}
Consequently, the set $\mathcal{M}_p$ is non-empty if and only if there exists $w\in\C^n\setminus\{0\}$ such that $\tsp\bar{w}Aw=\Re p_0-\tsp\bar{p}_\alpha Ap_\alpha$. This gives the three cases.

Set $x_0=\Re p_0-\tsp\bar{p}_\alpha Ap_\alpha$. According to the previous proposition, it suffices to consider the map
$$\left\{\!(y_0,v,w,a)\in\R\!\times\!\C^n\!\times\!(\C^n\!\setminus\!\{0\})\!\times\!\D/\ v=p_\alpha,y_0=\Im p_0,\frac{\tsp\bar{w}Aw}{1-|a|^2}=x_0\right\}\overset{\psi}{\to}\R\times\C^n$$
defined by
$$\psi(v,w,a,y_0)\!=\!\!(\Im h_0(1),h_\alpha(1))\!=\!\!\left(\!2\,\Im\left[\frac{\tsp\bar{p}_\alpha Aw}{1-a}\right]+2x_0\frac{\tsp\bar{w}Aw}{1-|a|^2}+\Im p_0,p_\alpha+\frac{w}{1-a}\right)\!.$$
Its differential at some point $(y_0,v,w,a)$ is defined on the tangent space
$$\left\{(0,0,w',a')\in\R\times\C^n\times\C^n\times\C\times/\ \frac{\tsp\bar{w'}Aw+\tsp\bar{w}Aw'}{1-|a|^2}+x_0\frac{\bar{a}a'+a\bar{a'}}{1-|a|^2}=0\right\}$$
by
$$\begin{array}{rl}
d\psi_{(y_0,v,w,a)}\!\!\!&\!\!\!(0,0,w',a')\\
 &=\left(\!2\,\Im\left[\frac{\tsp\bar{p}_\alpha Aw'}{1-a}+\frac{\tsp\bar{p}_\alpha Aw}{(1-a)^2}a'+\frac{x_0}{(1-a)^2}a'\right],\,\frac{w'}{1-a}+\frac{wa'}{(1-a)^2}\!\right).
\end{array}$$
Since $\mathcal{M}_p$ is non-empty, it is a submanifold of (real) dimension $2n+1$, and hence $d\psi_{(y_0,v,w,a)}$ is an isomorphism on $\R\times\C^n$ if and only if it is injective. Moreover,
\begin{eqnarray*}
d\psi_{(y_0,v,w,a)}(0,0,w',a')=0&\Leftrightarrow&\left\{\begin{array}{l}
 \frac{\tsp\bar{w}'Aw+\tsp\bar{w}Aw'}{1-|a|^2}+x_0\frac{\bar{a}a'+a\bar{a}'}{1-|a|^2}=0\\
 \Im\left[\frac{\tsp\bar{p}_\alpha Aw'}{1-a}+\frac{\tsp\bar{p}_\alpha Aw}{(1-a)^2a'}+\frac{x_0}{(1-a)^2}a'\right]=0\\
\frac{w'}{1-a}+\frac{w}{(1-a)^2}a'=0
\end{array}\right.\\
 &\Leftrightarrow&\left\{\begin{array}{l}
w'=-\frac{a'}{1-a}w\\
(-\frac{\bar{a}'}{1-\bar{a}}-\frac{a'}{1-a})\frac{\tsp\bar{w}Aw}{1-|a|^2}+x_0\frac{\bar{a}a'+a\bar{a}'}{1-|a|^2}=0\\
\Im(\frac{x_0}{(1-a)^2}a')=0
\end{array}\right.\\
 &\Leftrightarrow&\left\{\begin{array}{l}
w'=-\frac{a'}{1-a}w\\
x_0\left[-\frac{\bar{a}'}{1-\bar{a}}-\frac{a'}{1-a}+\frac{\bar{a}a'+a\bar{a}'}{1-|a|^2}\right]=0\\
x_0\,\Im\frac{a'}{(1-a)^2}=0
\end{array}\right.\,.
\end{eqnarray*}
In particular, if $x_0=0$ (that is, if $p\in Q$), $d\psi_{(v,w,a,y_0)}$ is not injective. If $x_0\not=0$, we can replace in the third row $\bar{a}'$ by $\frac{(1-\bar{a})^2}{(1-a)^2}$. Then $a'=0$ and $w'=0$, which concludes the proof.  
\end{proof}

\subsection{Stationary discs glued to a small perturbation of $Q$}\label{subs:StationaryDiscs}
\subsubsection{Discs glued to a small perturbation of $\PP(N^*Q)$}\label{subsubs:Glob}
The method consists in using a theorem of J. Globevnik \cite{Glob}:
 given a submanifold $E$ and a disc $f$ glued to $E$, and under the conditions that some integers depending on $E$ and $f$ are non-negative, the holomorphic discs near $f$ glued to a submanifold near $E$ form a $\kappa$-parameter family, where $\kappa$ is the Maslov index of $E$ along $f$. 
Since every stationary disc $h^0$ glued to a hypersurface $Q$ has a holomorphic lift $\hat{h}^0$ glued to the projectivization $\PP(N^*Q)$ of the conormal bundle of $Q$, we are going to apply the previous criterion to $E=\PP(N^*Q)$ and $\hat{h}^0$. Let us precise Globevnik's statement.\\

Let $f:\bd\to\C^N$ be the restriction to $\bd$ of a holomorphic disc, and assume that $f$ is of class $\mathcal{C}^\alpha$ and $\B\subset\R^{2N}$ is an open ball centered at 0. We suppose that there exist some functions $r_1,\hdots, r_N$ in $\mathcal{C}^2(\B)$ such that $dr_1\wedge\hdots\wedge dr_N$ does not vanish on $\B$, and 
$$M=\{\omega\in f(\zeta)+\B/\ \forall 1\le j\le N,\ r_j(\omega)=0\}$$
verifies $f(\zeta)\in M$ for all $\zeta\in\partial\D$. We also assume $T(\zeta)=T_{f(\zeta)}M$ to be totally real for all $\zeta\in\bd$.\\
Pick $\zeta\in\bd$, and denote by $G(\zeta)$ the invertible matrix $\displaystyle \left(\frac{\partial r_i}{\partial\bar{z}_j}(f(\zeta))\right)_{i,j}$. For any matrix $A(\zeta)$ whose columns generate $T(\zeta)$, each row of $G(\zeta)$ is orthogonal to each column of $A(\zeta)$:
$$\Re(\overline{G(\zeta)}A(\zeta)=0\Longleftrightarrow G(\zeta)\overline{A(\zeta)}=-\overline{G(\zeta)}A(\zeta)\Longrightarrow A(\zeta)\overline{A(\zeta)}^{-1}=-\overline{G(\zeta)}^{-1}G(\zeta).$$Set $B(\zeta)=A(\zeta)\overline{A(\zeta)}^{-1}$ for all $\zeta\in\bd$. The matrix $B$ does not depend on $A$. 

Moreover, one can find a Birkhoff factorization of $B$ \cite{Birkhoff}, that is some continuous matrix functions $B^+:\bar{\D}\to GL_N(\C)$ and $B^-:(\C\cup\{\infty\})\setminus\D\to GL_N(\C)$  such that 
$$\forall \zeta\in\bd,\ B(\zeta)=B^+(\zeta)\Lambda(\zeta)B^-(\zeta)$$
where $B^+$ and $B^-$ are holomorphic, and  
$\Lambda(\zeta)=\mathrm{diag}(\zeta^{\kappa_1},\hdots,\zeta^{\kappa_N})$.
The integers $\kappa_1\ge\hdots\ge\kappa_N$ do not depend on this factorization. They are called the {\em partial indices} of $B$ (see \cite{Vekua,CG} for more details).
We call the {\em partial indices of $M$ along $f$} the integers $\kappa_1\ge\hdots\ge\kappa_N$. Notice that they only depend on the bundle $\{T(\zeta)/\ \zeta\in\bd\}$. The {\em Maslov index}, or {\em total index}, of $M$ along $f$ is $\kappa=\sum_{1}^N\kappa_j$.

We only recall a slightly simpler version than the general result obtained by J. Globevnik.

\begin{theorem}{\bf \cite{Glob}}\label{theo:Globev} 
We assume that the previous conditions hold. For every $\rho=(\rho_1,\hdots,\rho_N)\in\mathcal{C}^2(\B)^N$ in a neighborhood of $r=(r_1,\hdots,r_N)$, we set
$$M_\rho=\{\omega\in f(\zeta)+\B/\ \forall 1\le j\le N,\ \rho_j(\omega)=0\}.$$
Assume that the partial indices of $M=M_r$ along $f$ are non-negative, and denote by $\kappa$ the Maslov index of $M$ along $f$. Then, there exist some open neighborhoods $V$ of $r$ in $\mathcal{C}^2(\B)^N$, $U$ of 0 in $\R^{\kappa+N}$, $W$ of $f$ in $\mathcal{C}^\alpha_{\C}(\bd)^N$, and a map $\mathcal{F}:V\times U\to\mathcal{C}^\alpha_{\C}(\bd)^N$ of class $\mathcal{C}^1$
such that
\begin{itemize}
\item $\mathcal{F}(r,0)=f$;
\item for all $(\rho,t)\in V\times U$, the map $\zeta\mapsto\mathcal{F}(\rho,t)(\zeta)-f(\zeta)$ is the boundary of a holomorphic disc glued to $M_\rho$;
\item there exists $\eta>0$ such that 
$\forall\rho\in V,\forall t_1,t_2\in U,\ ||\mathcal{F}(\rho,t_1)-\mathcal{F}(\rho,t_2)||\ge\eta|t_1-t_2|$,
in particular $\mathcal{F}(\rho,t_1)\not=\mathcal{F}(\rho,t_2)$ for $t_1\not= t_2$;
\item if $\dot{f}\in W$ is the boundary of a holomorphic disc glued to $M_\rho$, then there exists $t\in U$ such that $\dot{f}=\mathcal{F}(\rho,t)$.
\end{itemize}
\end{theorem}

When the partial indices are only supposed to be larger than -1, the result has been extended by Y.-G. Oh \cite{Oh95}.\\ 

Let $h$ be a non-constant stationary disc centered at $p=(p_0,0,\hdots,0)\notin Q$ and glued to $Q$. According to (\ref{CNScentre}), $w=h'_\alpha(0)\not=0$, which allows us to assume $(\tsp\bar{w}A)_n\not=0$. We denote by $\PP(N^*Q)$ the projectivization with respect to the $n-th$ coordinate.
Let $h^*$ be a regular lift of $h$, and 
$$\hat{h}^*=\left(\frac{h_0^*}{h_n^*},\hdots,\frac{h_{n-1}^*}{h_n^*}\right)$$ 
be its projectivization. Then $f=(h,\hat{h}^*)_{|\partial\Delta}$ is the boundary of a disc glued to $\PP(N^*Q)$, and 
\begin{equation}\label{exprf}
f(\zeta)\!=\!\!\left(\!\frac{\tsp{w}Aw}{1-|a|^2}\frac{1+a\zeta}{1-a\zeta}+iy_0,\frac{\zeta w_1}{1-a\zeta},\hdots,\frac{\zeta w_n}{1-a\zeta},\frac{\bar{a}-\zeta}{2(\tsp\bar{w}A)_n},\frac{(\tsp\bar{w}A)_1}{(\tsp\bar{w}A)_n},\hdots,\frac{(\tsp\bar{w}A)_{n-1}}{(\tsp\bar{w}A)_n}\!\right)\,
\end{equation}
in view of Remark \ref{exprcentre}.
In particular, $f$ is of class $\mathcal{C}^\infty$ up to the boundary, and does not depend on the choice of the regular lift $h^*$.

We want to apply the theorem of Globevnik to $M=\PP(N^*Q)$ along the holomorphic disc $f$. Note that $T(\zeta)=T_{f(\zeta)}(\PP(N^*Q))$ is totally real in view of Proposition \ref{prop:conormal}.
The difficulty is to compute the partial indices and the Maslov index of $\PP(N^*Q)$ along $f$. For the convenience of the reader, we prove the two following lemmas in the Appendices.

\begin{lemma}\label{indicespartiels}
The partial indices of $\PP(N^*Q)$ along $f$ are non-negative. 
\end{lemma}

\begin{lemma}\label{lm:indicetot}
The Maslov index $\kappa=\sum_{j=0}^{2n}\kappa_j$ of $\PP(N^*Q)$ along $f$ is equal to $2n+2$.
\end{lemma}

Hence we get: 

\begin{prop}\label{disquesproj}
Theorem \ref{theo:Globev} applies to our situation, and gives that the set of the discs we consider forms a real $[(2n+1)+\kappa]$-parameter family.
\end{prop}

\subsubsection{Discs glued to a small perturbation of $Q$}\label{subsubs:Perturbation}
Proposition \ref{disquesproj} shows that the set of all holomorphic discs near $f=(h,\hat{h}^*)$, glued to a small perturbation of $\PP(N^*Q)$ forms a manifold of (real) dimension $4n+3$. This gives us the description of the set of all holomorphic discs near $h$ glued to a small perturbation of $Q$ by means of the following lemma:  

\begin{lemma}{\bf(see \cite{L1,ST})}
Given $t$ sufficiently small and $\rho$ near $r$ with respect to the $\mathcal{C}^3$-topology, $\mathcal{F}(\rho,t)$ is the projectivization of a regular lift of some stationary disc.
\end{lemma}

\begin{proof}
We denote by $s$ the parameter $(\rho-r,t)$. Then, the disc $\mathcal{F}(\rho,t)$ is written in coordinates $(h_0^s,\hdots,h_n^s,H_0^s,\hdots,H_{n-1}^s)$. The holomorphic disc $h^s$ is glued to $M^s=\pi(\mathcal{M}^s)$, where $\mathcal{M}^s=\{\rho=0\}$ and $\pi$ is the canonical projection on the $n+1$ first coordinates. Hence $M^s=\{\rho_0=0\}$ and we may set
$$\phi^s:\bd\ni\zeta\mapsto\zeta\frac{\partial\rho_0}{\partial z_n}(h^s(\zeta)).$$
 For $s=0$, we obtain the initial disc $f$ and
for all $\zeta\in\bd$, $\phi^0(\zeta)=\frac{-1}{1-\bar{a}\bar{\zeta}}(\tsp\bar{w}A)_n$.
Note that the function $\bd\ni\zeta\mapsto\frac{-1}{1-\bar{a}\zeta}$ takes its values in the half-plane $\{\Re(z)<0\}$. Thus, for every sufficiently small $s$, the function $\phi^s$ also takes its values in some open half-plane which does not contains 0. This allows us to define 
$$\psi^s:\zeta\mapsto\log(\phi^s(\zeta))=\int_{[z_0;\phi^s(\zeta)]}\frac{\mathrm{d}z}{z}.$$
The function $\psi^s$ is of class $\mathcal{C}^\alpha$. Set $U^s=-\mathcal{T}(\Im\psi^s)$ and $\lambda^s=\mathrm{exp}\,(U^s-\Re\psi^s)$, where $\mathcal{T}$ is the Hilbert transform. The function $\lambda^s$ is positive valued, and we have on $\bd$:
$$\mathrm{log}\,(\lambda^s\phi^s)=(U^s-\Re\psi^s)+(\Re\psi^s+i\Im\psi^s)$$
which extends holomorphically to $\D$ in some function always denoted by $U^s+i\Im\psi^s$. Hence, $\mathrm{exp}\,(U^s+i\Im\psi^s)$ is a holomorphic extension of $\lambda^s\phi^s$ Ã  $\D$ that does not vanish in $\bar{\D}$. For all $0\le j\le n-1$,
$$\lambda^s\times\zeta\frac{\partial\rho_0}{\partial z_j}\circ h^s=\left(\lambda^s\times\zeta\frac{\partial\rho_0}{\partial z_n}\circ h^s\right)\times H_j^s$$
is a product of functions which extend holomorphically to $\D$. The function defined by $h^{*\,s}:=\lambda^s\frac{\partial\rho}{\partial z}\circ h^s$ is thus a regular lift of $h^s$. Hence, $(h^s,H^s)$ is the projectivization of $(h^s,h^{*\,s})$. Moreover, $\zeta h_n^{*\,s}$ does not vanish in $\bar{\D}$.
\end{proof} 

This proves that the map $(h,h^*)\mapsto (h,\hat{h}^*)$ is onto. 
For any stationary disc $h$ glued to $M$, and $h^{*\,(1)}$, $h^{*\,(2)}$ two regular lifts of $h$: 
$$\forall\zeta\in\bd,\ \exists\lambda(\zeta)\in\R^*/\ h^{*\,(2)}(\zeta)=\lambda(\zeta)h^{*\,(1)}(\zeta).$$
Since $\zeta h^{*\,(j)}_n$ is holomorphic and does not vanish in $\bar{\D}$ by hypothesis, the function $\frac{\zeta h_n^{*\,(2)}}{\zeta h_n^{*\,(1)}}=\lambda$ extends holomorphically to $\D$: hence $\lambda$ is a real non-zero constant. Then the stationary disc $h$ has an only regular lift up to multiplication by a real non-zero constant, and two regular lifts have the same image under $(h,h^*)\mapsto (h,\hat{h}^*)$. Consequently, the map $h\mapsto (h,\hat{h}^*)$ is a one-to-one correspondence. Finally, we get that the set of all the stationary discs near $h^0$ glued to a small deformation of $Q$ form a manifold of (real) dimension $4n+3$. More precisely, we have proved:

\begin{theorem}\label{theo:parametrisation}
Let $\B\subset\R^{2n}$ be an open ball in a neighborhood of 0, $A$ an $n$-sized non-degenerate hermitian matrix, and $Q=\{0=r(z)=\Re z_0-\tsp\bar{z_\alpha}Az_\alpha\}$. Let $h^0$ be a stationary disc glued to $Q$ such that $h^0(0)=(p_0,0,\hdots,0)\notin Q$.
Then there exist some open neighborhoods $V$ of $r$ in $\mathcal{C}^3(\B,\R)$, $U$ of 0 in $\R^{4n+3}$, $W$ of $h^0$ in $\mathcal{C}^\alpha_{\ccc}(\bd)^{n+1}$, and a map $\mathcal{H}:V\times U\to\mathcal{C}^\alpha_{\ccc}(\bd)^{n+1}$ of class $\mathcal{C}^1$ such that:
\begin{itemize}
\item $\mathcal{H}(r,0)=h^0$;
\item for all $(\rho,t)\in V\times U$, $\mathcal{H}(\rho,t)$ is the boundary of a stationary disc glued to $M$, where $M=\{\rho=0\}$;
\item if $t_1\not= t_2$, $\mathcal{H}(\rho,t_1)\not=\mathcal{H}(\rho,t_2)$;
\item if $\dot{h}\in W$ is the boundary of a stationary disc glued to $M=\{\rho=0\}$, with $\rho\in V$, then there exists $t\in U$ such that $\dot{h}=\mathcal{H}(\rho,t)$.
\end{itemize}
\end{theorem}

Note that, by Corollary \ref{CNSexistence}, there exists a non-constant stationary disc $h^0$ glued to $Q$ and such that $h^0(0)=(p_0,0,\hdots,0)$ if and only if the following condition holds:\\

\noindent\textit{{\bf Condition \boldmath $\ast$\unboldmath} The point $p=(p_0,0,\hdots,0)$ verifies $p\notin Q$. Moreover, $\Re p_0>0$ (resp. $\Re p_0<0$) if $A$ is positive definite (resp. negative definite).}\\

\begin{cor}\label{passageQM}
For $\rho$ near $r$ with respect to the $\mathcal{C}^3$-topology, the map $\Phi^{-1}$ defined in (\ref{parametrisationexplicite}) is always a local $\mathcal{C}^1$-diffeomorphism from $\R\times\C^n\times(\C^n\setminus\{0\})\times\D$ to the set $\mathcal{M}^\rho$ of all non-constant stationary discs glued to $M=\{\rho=0\}$.
\end{cor}

\begin{proof}
We use the notations of Theorem \ref{theo:parametrisation} and Proposition \ref{forme}. For all $\rho\in V$, we define the map 
$$\begin{array}{cccl}
\Theta_\rho:&U&\to&\R\times\C^n\times(\C^n\setminus\{0\})\times\D\\
            &t&\mapsto&(y_0^{(\rho,t)},v^{(\rho,t)},w^{(\rho,t)},a^{(\rho,t)})
\end{array}
$$
where $y_0^{(\rho,t)},v^{(\rho,t)},w^{(\rho,t)},a^{(\rho,t)}$ are the parameters defining the disc $\mathcal{H}(\rho,t)$. We can assume that any disc $h=\mathcal{H}(\rho,t)$ verifies $(\tsp\bar{h}'_\alpha(0)A)_n\not=0$. The map $\Theta_\rho$ is of class $\mathcal{C}^1$. Moreover, the map $\mu:(\rho,t)\mapsto d{(\Theta_\rho)}_t$ is continuous from $V\times U$ to the Banach space $\mathcal{L}_c(\R^{4n+3})$ of continuous linear maps. Since $\mu(r,0)$ is invertible, we can assume that the differential of $\Theta_\rho$ is invertible at any point. This gives the conclusion by means of the inverse function theorem.
\end{proof}

\subsection{Discs centered at $p$}\label{subs:FixedCenter}
We can now give two parametrizations of the discs centered at some fixed point and glued to a small perturbation of $Q$. We begin with studying discs glued to $Q$.

\begin{prop}
Pick $p=(p_0,0,\hdots,0)$ verifying Condition $\ast$.
Then the map that associates to any non-constant stationary disc $h$, glued to $Q$ and centered at $p$, the point $h(1)$ (resp. the vector $h'(0)$), is a local diffeomorphism on its image.
\end{prop}

\begin{proof}
The map $h\mapsto h(1)$ has been studied yet in Corollary \ref{CNSexistence}. For $h\mapsto h'(0)$, according to Proposition \ref{forme} and Remark \ref{exprcentre}, it suffices to consider the map 
$$\Psi:(\C^n\setminus\{0\})\times\D\ni (w,a)\mapsto \left(2a\frac{\tsp\bar{w}Aw}{1-|a|^2},w\right).$$
If $\Psi(w,a)=\Psi(\dot{w},\dot{a})$, then $w=\dot{w}$ and
$$\frac{a}{1-|a|^2}=\frac{\dot{a}}{1-|\dot{a}|^2}\Leftrightarrow\left\{\begin{array}{l}\frac{|a|}{1-|a|^2}=\frac{|\dot{a}|}{1-|\dot{a}|^2}\\\mathrm{Arg}\,a=\mathrm{Arg}\,\dot{a}\end{array}\right.\Leftrightarrow\left\{\begin{array}{l}|a|=|\dot{a}|\\\mathrm{Arg}\,a=\mathrm{Arg}\,\dot{a}.\end{array}\right.$$
Hence $\Psi$ is injective and induces a homeomorphism on its image. Moreover, 
$$d\Psi_{(w,a)}(w',a')=\left(\begin{array}{ccc}2a\frac{\tsp\bar{w}'Aw+\tsp\bar{w}Aw'}{1-|a|^2}&+&2\frac{\tsp\bar{w}Aw}{(1-|a|^2)^2}\,(a'+a^2\bar{a}')\\ w'&+&0\end{array}\right),$$
thus $d\Psi$ is injective at every point.
\end{proof}

By means of arguments similar to those of the proof of Corollary \ref{passageQM}, one gets Theorem \ref{2diffeo}.

\begin{remark}\label{rem:prec}
 More precisely, the map $h\mapsto(\Im h_0(1),h_\alpha(1))$ is a local diffeomorphism on $\R\times\C^n$. 
\end{remark}

\section{Uniqueness property}\label{s:U}
Let $M$ be a small perturbation of the hyperquadric $Q$. As soon as $z\in Q$ is met by some stationary disc, Theorem \ref{2diffeo} provides, in a neighbourhood of $z$, a foliation of $M$ by the boundaries of stationary discs. The following lemma determines such points $z$:

\begin{lemma}
Assume that the conditions of Theorem \ref{2diffeo} hold. We denote by $\delta^\rho$ the map defined on $\mathcal{M}_p^{(\rho)}$ by $\delta^\rho(h)=h(1)$.\\
Then the image of the map $\delta^r$ is exactly $\{z\in Q/\ \Re z_0\Re p_0>0\}$. In particular, if $A$ is positive definite or negative definite, then every $z\in Q\setminus\{0\}$ may be written under the form $z=h(1)$ for some $h\in\mathcal{M}_p^{(r)}$.
\end{lemma}

\begin{proof}
Let $h$ be a stationary disc glued to $Q$ and centered at 
$p\!=\!\!\left(\frac{\tsp\bar{w}Aw}{1-|a|^2}+iy_0,0\right)$.
Thus, for $z\in Q$:
$$z\!=\!h(1)\!=\!\!\left(\frac{\tsp\bar{w}Aw}{1-|a|^2}\,\frac{1+a}{1-a}+iy_0,w\frac{1}{1-a}\right)\Leftrightarrow\left\{\!\!\begin{array}{l} w=(1-a)z_\alpha\\ z_0\not=-\Re p_0+i\,\Im p_0\ \mathrm{and}\ a=\frac{z_0-p_0}{z_0+\bar{p}_0}\,. \end{array}\right.$$
It admits a solution if and only if $z_0\not=-\bar{p}_0$ and
$$\left|\frac{z_0-p_0}{z_0+\bar{p}_0}\right|<1\Leftrightarrow \frac{(\Re z_0-\Re p_0)^2+(\Im z_0-\Im p_0)^2}{(\Re z_0+\Re p_0)^2+(\Im z_0-\Im p_0)^2}<1\Leftrightarrow \Re z_0\Re p_0>0,$$
which concludes the proof.
\end{proof}

\noindent In view of Remark \ref{rem:prec}, this proves that every point of $M$ near such a point $z\in Q$ has a neighbourhood which is foliated by the boundaries of stationary discs.\\

Theorem \ref{2diffeo} allows us to construct a local analogue of the Kobayashi indicatrix, considering the real hypersurface $I(M)=\{h'(0)/\ h\in\mathcal{M}_p^{(\rho)}\}$ (we always assume that the conditions of the theorem hold).

\begin{figure}
\begin{picture}(12,8.5)
\qbezier(0,4.75)(0.6,5.25)(2.5,5.2)
\qbezier(2.5,5.2)(4.3,5.3)(4.6,5.8)
\qbezier(1,6.55)(1.6,7.05)(3.5,7)
\qbezier(3.5,7)(5.3,7.1)(5.6,7.6)
\qbezier(0,4.75)(0,4.95)(1,6.55)
\qbezier(4.6,5.8)(4.6,6)(5.6,7.6)

\put(2.8,6.2){\circle*{0.07}}
\put(2.8,8){\circle*{0.07}}

\qbezier(6.9,5.1)(7.7,5.5)(9.6,5.3)
\qbezier(9.6,5.3)(11.6,5)(11.9,5.3)
\qbezier(7.3,7)(7.9,7.3)(9.8,7.2)
\qbezier(9.8,7.2)(11.8,6.9)(12.2,7.1)
\qbezier(6.9,5.1)(6.9,5.1)(7.3,7)
\qbezier(11.9,5.3)(11.9,5.3)(12.2,7.1)

\put(9.7,6.27){\circle*{0.07}}
\put(9.7,8){\circle*{0.07}}

\qbezier[30](2.6,6.45)(2,6.2)(2.5,5.9)
\qbezier[10](2.6,6.45)(2.8,6.5)(3,6.45)
\qbezier[30](3,6.45)(3.6,6.2)(3,5.9)
\qbezier[10](2.5,5.9)(2.8,5.8)(3,5.9)
\qbezier[80](2.28,6.2)(2.8,9.8)(3.3,6.2)

\qbezier(2.2,6.65)(1.5,6.2)(2,5.7)
\qbezier(2.2,6.65)(2.8,6.9)(3.4,6.65)
\qbezier(3.4,6.65)(4.1,6.2)(3.6,5.7)
\qbezier(2,5.7)(2.8,5.2)(3.6,5.7)
\qbezier(1.8,6.2)(2.8,9.8)(3.8,6.2)
\put(2,5.7){\circle*{0.07}}
\put(1.4,5.5){$_{h(1)}$}

\put(0.5,5.2){$M$}
\put(2.8,8.2){$p$}
\put(3,5.9){\circle*{0.07}}
\put(3.05,5.8){$_z$}

\qbezier[30](9.5,6.5)(8.9,6.3)(9.4,6)
\qbezier[10](9.5,6.5)(9.7,6.55)(9.9,6.5)
\qbezier[30](9.9,6.5)(10.5,6.3)(9.95,6)
\qbezier[10](9.4,6)(9.7,5.9)(9.95,6)
\qbezier[80](9.18,6.3)(9.7,9.7)(10.21,6.3)

\qbezier(9.35,6.62)(8.5,6.3)(9.2,5.9)
\qbezier(9.35,6.62)(9.7,6.72)(10.05,6.62)
\qbezier(10.05,6.62)(10.9,6.3)(10.15,5.9)
\qbezier(9.2,5.9)(9.7,5.7)(10.15,5.9)
\qbezier(8.9,6.3)(9.7,9.7)(10.5,6.3)
\put(10.39,6.05){\circle*{0.07}}
\put(10.5,5.95){$_{F\circ h(1)}$}

\put(7.2,5.45){$M'$}
\put(9.65,8.15){$F(p)$}

\put(5.65,6.85){\vector(1,0){1}}
\put(6,7){$F$}

\put(2.7,4.4){\vector(0,-1){1.5}}
\put(2.9,3.6){$\begin{array}{c}h(1)\\ \downarrow \\ h'(0)\end{array}$}
\put(3.4,3.89){$_{\!_{-}}$}

\put(9.7,4.4){\vector(0,-1){1.5}}
\put(9.9,3.6){$\begin{array}{c}h(1)\\ \downarrow \\ h'(0)\end{array}$}
\put(10.4,3.89){$_{\!_{-}}$}

\put(5.55,1.7){\vector(1,0){1}}
\put(5.7,1.87){$dF_p$}

\put(2.7,1.7){\circle{2}}
\qbezier(2,1.7)(2.7,1)(3.4,1.7)
\qbezier[20](2,1.7)(2.7,2.3)(3.4,1.7)
\put(2.3,0.4){$I(M)$}

\put(9.7,1.7){\circle{2}}
\qbezier(9,1.7)(9.7,1)(10.4,1.7)
\qbezier[20](9,1.7)(9.7,2.3)(10.4,1.7)
\put(9.3,0.4){$I(M')$}

\end{picture}
\caption{}
\label{dessin}
\end{figure}

The diffeomorphism $h'(0)\mapsto h(1)$ is well-defined, and Figure \ref{dessin} shows that it commutes with the biholomorphisms in the following sense. Of course the statement will be local. 
 If $p=(p_0,0)\notin Q$ and $F(p)=(p'_0,0)\notin Q'$ verify Condition $\ast$, the following diagram commutes:

$$\begin{array}{rccc}
F:&(\Omega,M)&\longrightarrow&(\Omega',M')\\
  &h(1)&\mapsto&F\circ h(1)\\
  &\updownarrow& \circlearrowleft&\updownarrow\\
  &h'(0)&\underset{dF_p}{\longrightarrow}&dF_p(h'(0)). 
\end{array}$$

\noindent Hence $F$ is uniquely determined by its differential at $p$. Let us summarize:

\begin{theorem}\label{theo:unicite}
Let $M$ and $M'$ be two real hypersurfaces of $\C^{n+1}$ given respectively by 
$$0=\Re z_0-\sum_{i,j\ge 1}a_{i,j}\bar{z}_iz_j+O(||(\Im z_0,z_\alpha)||^3)$$
and
$$0=\Re z_0-\sum_{i,j\ge 1}a'_{i,j}\bar{z}_iz_j+O(||(\Im z_0,z_\alpha)||^3),$$
where $A=(a_{i,j})$ and $A'=(a'_{i,j})$ are non-degenerate hermitian matrices. We denote by $Q$ and $Q'$ the hyperquadrics associated to $A$ and $A'$.\\
Let $F:\Omega\to\Omega'$ be a biholomorphism such that $F(M)\subset M'$. We assume that there exists $p=(p_0,0)\notin Q$ such that $F(p)=(p'_0,0)\notin Q'$ and $p,\ F(p)$ verify Condition $\ast$. Finally, pick $z\in M$ such that $\Re z_0\Re p_0>0$.\\
If $M$ and $M'$ are sufficiently near $Q$ and $Q'$, then $F$ is determined in a neighborhood of $z$ by its differential at $p$.
\end{theorem}

The neighborhood of the equations of $Q$ and $Q'$ is understood in the $\mathcal{C}^3$-sense, and depends on $p$ and $z$. The hypotheses on $F$ just assure that Theorem \ref{2diffeo} applies to the source and target spaces.

\appendix
\section{Proof of Lemma \ref{indicespartiels}}\label{s:IndicesPart}
Here we use the notations of Section \ref{subsubs:Glob}. We begin with determining the equations of $\PP(N^*Q)$. Since the fibers of $N^*Q$ are generated by $\partial r$, a point $(z,t)\in\C^{n+1}\times\C^n$ is in $\PP(N^*Q)$ if and only if 
\begin{eqnarray*}
0&=&r(z)\\
0&=&\frac{\partial r}{\partial z_n}(z)t_j-\frac{\partial r}{\partial z_j}(z)\quad\mathrm{for\ all}\ j=0,\hdots, n-1\,.
\end{eqnarray*}

By separating the real and imaginary parts, we obtain $2n+1$ equations $r_0=\hdots=r_{2n}=0$ (with $r_0=r$), and the matrix $G(\zeta)=\left(\frac{\partial r_i}{\partial\bar{z}_j}(f(\zeta))\right)_{i,j}$ becomes
$$\left(\!\!\!\begin{array}{cccccccc}
1/2\!\!&\!\! -L_1z_\alpha\!\!&\!\!\hdots\!\!&\!\!-L_nz_\alpha\!\!&\!0\!&\!\!0\!\!&\!\!\hdots\!\!&\!\! \\
0\!\!&\!\!-\bar{A}_{n,1}t_1+\bar{A}_{1,1}\!\!&\!\!\hdots\!\!&\!\!-\bar{A}_{n,n}t_1+\bar{A}_{1,n}\!\!&\!\!0\!\!&\!\!-L_nz_\alpha\!\!&\!\! \!\!&\!\! \\
\vdots\!\!&\!\!\vdots\!\!&\!\! \!\!&\!\!\vdots\!\!&\!\!\vdots\!\!&\!\! \!\!&\!\!\ddots\!\!&\!\! \\
0\!\!&\!\!-\bar{A}_{n,1}t_{n-1}+\bar{A}_{n-1,1}\!\!&\!\!\hdots\!\!&\!\!-\bar{A}_{n,n}t_{n-1}+\bar{A}_{n-1,n}\!\!&\!\!0\!\!&\!\! \!\!&\!\! \!\!&\!\!-L_nz_\alpha\\
0\!\!&\!\!-\bar{A}_{n,1}t_0\!\!&\!\!\hdots\!\!&\!\!-\bar{A}_{n,n}t_0\!\!&\!\!-L_nz_\alpha\!\!&\!\!0\!\!&\!\!\hdots\!\!&\!\!0 \\
0\!\!&\!\!i(-\bar{A}_{n,1}t_1+\bar{A}_{1,1})\!\!&\!\!\hdots\!\!&\!\!i(-\bar{A}_{n,n}t_1+\bar{A}_{1,n})\!\!&\!\!0\!\!&\!\!iL_nz_\alpha\!\!&\!\! \!\!&\!\!  \\
\vdots\!\!&\!\!\vdots\!\!&\!\! \!\!&\!\!\vdots\!\!&\!\!\vdots\!\!&\!\! \!\!&\!\!\ddots\!\!&\!\! \\
0\!\!&\!\!\!\!i(-\bar{A}_{n,1}t_{n-1}+\bar{A}_{n-1,1})\!\!&\!\!\hdots\!\!&\!\!i(-\bar{A}_{n,n}t_{n-1}+\bar{A}_{n-1,n})\!\!&\!\!0\!\!&\!\! \!\!&\!\! \!\!&\!\!iL_nz_\alpha\\
0\!\!&\!\!-i\bar{A}_{n,1}t_0\!\!&\!\!\hdots\!\!&\!\!-i\bar{A}_{n,n}t_0\!\!&\!\!iL_nz_\alpha\!\!&\!\!0\!\!&\!\! \hdots\!\!&\!\!0
\end{array}\!\!\!\!\!\right)
$$
where $L_j$ denotes the $j-th$ row of the matrix $A$. Let us notice that
$$\left(\!\!\begin{array}{ccc}
(-\bar{A}_{n,1}t_1+\bar{A}_{1,1})\!&\!\hdots\!&\!(-\bar{A}_{n,n}t_1+\bar{A}_{1,n})\\
\vdots\!&\! \!&\!\vdots\\
(-\bar{A}_{n,1}t_{n-1}+\bar{A}_{n-1,1})\!&\!\hdots\!&\!(-\bar{A}_{n,n}t_{n-1}+\bar{A}_{n-1,n})\\
-\bar{A}_{n,1}t_0\!&\!\hdots\!&\!-\bar{A}_{n,n}t_0
\end{array}\!\!\right)\!=\!\left(\!\!\begin{array}{cccc}
1 \!&\! \!&\! \!&\!-t_1\\
 \!&\! \ddots\!&\! \!&\!\vdots\\
\!&\! \!&\! 1\!&\!-t_{n-1}\\
0\!&\!\hdots\!&\!0\!&\!-t_0
\end{array}\!\!\right)\!\!\times\!\bar{A}.$$
The right multiplication by the constant matrix $\left(\begin{array}{ccc}1&0&0 \\0&\bar{A}^{-1}&0\\ 0&0&I_n  \end{array}\right)$ does not change the partial indices, and gives us the matrix
\begin{equation}\label{mat}
\left(\begin{array}{ccccccccc}
1/2&-z_1&\hdots&-z_{n-1}&-z_n&0& &\hdots&0\\
0&1& & &-t_1&0&-L_nz_\alpha& & \\
\vdots& &\ddots& &\vdots&\vdots& &\ddots& \\
0& & &1&-t_{n-1}&0& & &-L_nz_\alpha\\
0&0&\hdots&0&-t_0&-L_nz_\alpha&0&\hdots&0\\
0&i& & &-it_1&0&iL_nz_\alpha& & \\
\vdots& &\ddots& &\vdots&\vdots& &\ddots& \\
0& & &i&-it_{n-1}&0& & &iL_nz_\alpha\\
0&0&\hdots&0&-it_0&iL_nz_\alpha&0&\hdots&0
\end{array}\right)\,.
\end{equation}
Replacing $f=(h,\hat{h}^*)$ by its expression, one gets, in view of (\ref{exprf}),
$z_\alpha=h_\alpha(\zeta)=\frac{\zeta}{1-a\zeta}w$, $t_0=-\frac{\zeta-\bar{a}}{2(\tsp\bar{w}A)_n}$ and for all $1\le j\le n-1$, $t_j=\frac{(\tsp\bar{w}A)_j}{(\tsp\bar{w}A)_n}$. Hence
$$L_nz_\alpha=\left(A\times\frac{\zeta}{1-a\zeta}w\right)_n=\frac{\zeta}{1-a\zeta}(Aw)_n=\frac{\zeta}{1-a\zeta}\tsp(Aw)_n=\frac{\zeta}{1-a\zeta}(\tsp w\bar{A})_n.$$
We number the columns of the matrix (\ref{mat}) from 0 to $2n$. By multiplying $C_0$ by 2 and $C_{n+1},\hdots, C_{2n}$ by $1/(\tsp w\bar{A})_n$, and then permuting the rows, we obtain:
$$\left(\begin{array}{ccccccccccc} 
1&-z_1&-z_2&\hdots&-z_{n-1}&-z_n&0&0&\hdots&\hdots&0\\
0&0&0&\hdots&0&-t_0&-\xi&0&\hdots&\hdots&0\\
0&0&0&\hdots&0&-it_0&i\xi&0&\hdots&\hdots&0\\
0&1&0&\hdots&0&-t_1&0&-\xi&\hdots&\hdots&0\\
0&i&0&\hdots&0&-it_1&0&i\xi&\hdots&\hdots&0\\
 & & &\ddots &\vdots&\vdots&\vdots & & &\ddots & \\
 &(0) & &\ddots &1&-t_{n-1}&0& & &\ddots  &-\xi\\
 & & & &i&-it_{n-1}&0&  &(0) & &i\xi
\end{array}\right)\,,$$
where $\xi=\frac{\zeta}{1-a\zeta}$.
Since each $t_j,\ 1\le j\le n-1$ is constant, the operations 
$C_1:=C_n+\sum_{j=1}^{n-1}t_jC_j$ $C_2:=C_{n+1}$
and $C_{2j+1}=C_j$, $C_{2j+2}=C_{n+j}$ if $j\ge 1$ do not change the partial indices and give
\begin{eqnarray}\label{matricecalcul}
\left(\begin{array}{cccccccc}
1&(-z_n-\sum_{j=1}^{n-1}t_jz_j)&0&-z_1&0&\hdots&-z_{n-1}&0\\
0&-t_0&-\xi& & & & & \\
0&-it_0&i\xi& & & & & \\
 & & &1&-\xi& &(0) & \\
 & & &i&i\xi& & & \\
 & & & &    &\ddots& & \\
 &(0) & & &    &      &1&-\xi\\
 & & & &    &      &i&i\xi
\end{array}\right)\,.
\end{eqnarray}
Along $f$, one has
$\displaystyle\sum_{j=1}^n\frac{\partial r}{\partial z_j}(f(\zeta))\times f_j(\zeta)=-\tsp\bar{h}_\alpha(\zeta)Ah_\alpha(\zeta)=-\Re(h_0(\zeta))$ and
$$-z_n-\sum_{j=1}^{n1}t_jz_j=-\frac{\left(\frac{\partial r}{\partial z_n}z_n+\sum_{j=1}^{n-1}\frac{\partial r}{\partial z_j}z_j\right)}{\frac{\partial r}{\partial z_n}}=\frac{\Re(h_0(\zeta))}{\frac{\partial r}{\partial z_n}}=2t_0\Re(h_0(\zeta)).$$ 
The $(3\times 3)$-sized upper left block in the matrix (\ref{matricecalcul}) is thus equal to 
$$\left(\begin{array}{ccc}1&2t_0\Re z_0&0\\ 0&-t_0&-\xi\\ 0&-it_0&i\xi\end{array}\right),$$ 
where $\xi=\frac{\zeta}{1-a\zeta}$ and $t_0$ is under the form $\gamma\zeta(1-\bar{a}\bar{\zeta})$ with $\gamma$ some non-zero constant. The right multiplication by 
$\left(\begin{array}{ccc}1&0&0\\0&\frac{1}{\gamma(1-\bar{a}\bar{\zeta})}&0\\0&0&I_{2n-1}\end{array}\right)$ 
allows us to consider, in place of the matrix (\ref{matricecalcul}), the matrix
$$G'(\zeta)\!=\!\!\left(\!\!\begin{array}{cccccccc}
1&\frac{\zeta}{|1-a\zeta|^2}\times 2\tsp\bar{w}Aw&0\!&\!\frac{-\zeta}{1-a\zeta}w_1\!&\!0\!&\!\hdots\!&\frac{-\zeta}{1-a\zeta}w_{n-1}\!&\!0\\
0&-\zeta&\frac{-\zeta}{1-a\zeta}\!&\! \!&\! \!&\! \!&\! \!&\! \\
0&-i\zeta&i\frac{\zeta}{1-a\zeta}& \!&\! \!&\! \!&\! \!&\! \\
 & & \!&\!1\!&\!\frac{-\zeta}{1-a\zeta}\!&\! \!&\!(0) \!&\! \\
 & & \!&\!i\!&\!i\frac{\zeta}{1-a\zeta}\!&\! \!&\! \!&\! \\
 & & \!&\! \!&\!  \!&\!\ddots\!&\! \!&\! \\
 &(0) & \!&\! \!&\! \!&\! \!&\!1\!&\!\frac{-\zeta}{1-a\zeta}\\
 & & \!&\! \!&\! \!&\! \!&\!i\!&\!i\frac{\zeta}{1-a\zeta}
\end{array}\!\!\right).$$ 

We write $G'(\zeta)=\left(\begin{array}{cc}\alpha&\beta\\ 0&\gamma\end{array}\right)$, where the entries are $\alpha\in\mathcal{M}_3(\C)$, $\beta\in\mathcal{M}_{3,2n-2}(\C)$ and $\gamma\in\mathcal{M}_{2n-2}(\C)$. Then
\begin{eqnarray*}
\lefteqn{\overline{G'(\zeta)^{-1}}G'(\zeta)\!\!=\!\!\left(\begin{array}{cc}\bar{\alpha}^{-1}\alpha\ &\ \bar{\alpha}^{-1}(\beta-\bar{\beta}\bar{\gamma}^{-1}\gamma)\\ 0\ &\ \bar{\gamma}^{-1}\gamma\end{array}\right)}\\
 \!\!&\!\!=\!\!&\!\!\left(\begin{array}{cccccc}
1&\frac{\zeta}{|1-a\zeta|^2}2\tsp\bar{w}Aw&\frac{-\zeta}{|1-a\zeta|^2(1-a\zeta)}2\tsp w\bar{A}\bar{w}&\frac{-\zeta}{1-a\zeta}w_1&\frac{-1}{|1-a\zeta|^2}\bar{w}_1&\hdots \\
0&0&\frac{\zeta^2}{1-a\zeta}& & &(0) \\
0&\zeta^2(1-\bar{a}\bar{\zeta})&0& & & \\
 & & &0&\frac{-\zeta}{1-a\zeta}&  \\
 & & &-\zeta(1-\bar{a}\bar{\zeta})&0&  \\
 &(0)& & & &\ddots 
\end{array}\right)\!.
\end{eqnarray*}
Let us number the rows and the columns of this matrix from $0$ to $2n$, and multiply $C_1$ by $(1-\bar{a}\bar{\zeta})$ and $L_1$ by $(1-a\zeta)$. We also multiply $C_j$ by $-1/(1-\bar{a}\bar{\zeta})$ and $L_j$ by $-(1-a\zeta)$ for all $j$ odd and larger than 3. Hence, we get $\overline{Q^{-1}}\overline{G'(\zeta)^{-1}}G'(\zeta)Q(\zeta)$, where $Q$ is the $(2n+1)$-sized diagonal matrix
$$Q=\mathrm{diag}\left(1,\frac{1}{1-\bar{a}\bar{\zeta}},1,\frac{-1}{1-\bar{a}\bar{\zeta}},1,\hdots,\frac{-1}{1-\bar{a}\bar{\zeta}},1\right).$$
Thus the partial indices of $\PP(N^*Q)$ along $f$ are those of the matrix $B(\zeta)=\overline{(G'(\zeta)Q(\zeta))}^{-1}(G'(\zeta)Q(\zeta))$, that is:
\begin{equation}\label{B}
B(\zeta)\!=\!\left(\begin{array}{cccccc}
1&\frac{\zeta}{|1-a\zeta|^2(1-\bar{a}\bar{\zeta})}2\tsp\bar{w}Aw&\frac{-\zeta}{|1-a\zeta|^2(1-a\zeta)}2\tsp w\bar{A}\bar{w}&\frac{\zeta w_1}{|1-a\zeta|^2}&\frac{-\bar{w}_1}{|1-a\zeta|^2}&\hdots\\
0&0&\zeta^2& &  \\
0&\zeta^2&0& & &(0)  \\
 & & &0&\zeta&  \\
 & & &\zeta&0&  \\
 &(0)& & & &\ddots \\
\end{array}\right).
\end{equation}
We use the following lemma:

\begin{lemma}{\bf (\cite{Glob}, Lemma 5.1)}\ 
Let $G'':\partial\Delta\to GL_{2n+1}(\C)$ of class $\mathcal{C}^\alpha$ ($0<\alpha<1$), and denote by $\kappa_0\ge\hdots\ge\kappa_{2n}$ the partial indices of the function $\zeta\mapsto G''(\zeta)\overline{G''(\zeta)}^{-1}$. Then there exists some function $\theta:\bar{\D}\to GL_n(\C)$ of class $\mathcal{C}^\alpha$, holomorphic on $\D$, such that
$$\forall\zeta\in\partial\D,\ G''(\zeta)\overline{G''(\zeta)}^{-1}=\theta(\zeta)\left(\begin{array}{ccc}\zeta^{\kappa_0}& &(0) \\ &\ddots& \\ (0)& &\zeta^{\kappa_{2n}}\end{array}\right)\overline{\theta(\zeta)}^{-1}.$$
\end{lemma}

\noindent Applying this result to the matrix $G''=\overline{(G'Q)}^{-1}$, we get some function $P$ from $\bar{\D}$ to $GL_{2n+1}(\C)$, holomorphic on $\D$, such that 
$$\forall\zeta\in\bd,\ P(\zeta)B(\zeta)=\left(\begin{array}{ccc}\zeta^{\kappa_0}& &(0) \\ &\ddots& \\ (0)& &\zeta^{\kappa_{2n}}\end{array}\right)\overline{P(\zeta)}.$$
In particular, if we denote by $l=(l_0,\hdots,l_{2n})$ the last row of $P$, we get that for all $\zeta\in\bd$, 
\begin{eqnarray*}
\lefteqn{l(\zeta)B(\zeta)=\zeta^{\kappa_{2n}}\overline{l(\zeta)}\ \ \ \ \ \ \ \ \ \ \ \ \ \ \ \ \ }\\
&\Longleftrightarrow& (S)\left\{\begin{array}{l}
l_0=\zeta^{\kappa_{2n}}\bar{l}_0\\
l_0\times\frac{2\tsp\bar{w}Aw\zeta}{|1-a\zeta|^2(1-\bar{a}\bar{\zeta})}+l_2\zeta^2=\zeta^{\kappa_{2n}}\bar{l}_1\\
l_0\times\frac{-2\tsp wA\bar{w}\zeta}{|1-a\zeta|^2(1-a\zeta)}+l_1\zeta^2=\zeta^{\kappa_{2n}}\bar{l}_2\\
\forall 1\le j\le n-1,\ \left\{\begin{array}{l} l_0\times\frac{\zeta w_j}{|1-a\zeta|^2}+\zeta l_{2j+2}=\zeta^{\kappa_{2n}}\bar{l}_{2j+1}\\
l_0\times\frac{-\bar{w}_j}{|1-a\zeta|^2}+\zeta l_{2j+1}=\zeta^{\kappa_{2n}}\bar{l}_{2j+2}.\end{array}\right.
\end{array}\right.
\end{eqnarray*}

If $l_0$ is not identically zero, the first row of this system implies $\kappa_{2n}\ge 0$ since $l_0$ is holomorphic and $\bar{l}_0$ is anti-holomorphic.

If $l_0\equiv 0$, there exists some $j\ge 1$ such that the function $l_j$ does not vanish identically. The system $(S)$ becomes
$$\left\{\begin{array}{l}
l_2\zeta^2=\zeta^{\kappa_{2n}}\bar{l}_1\quad \mathrm{and}\quad l_1\zeta^2=\zeta^{\kappa_{2n}}\bar{l}_2\\
\zeta l_4=\zeta^{\kappa_{2n}}\bar{l}_3\quad \mathrm{and}\quad \zeta l_3=\zeta^{\kappa_{2n}}\bar{l}_4\\
\vdots\\
\zeta l_{2n}=\zeta^{\kappa_{2n}}\bar{l}_{2n-1}\quad \mathrm{and}\quad \zeta l_{2n-1}=\zeta^{\kappa_{2n}}\bar{l}_{2n}.
\end{array}\right.$$
Consequently, reasoning for $l_j$ as previously, we obtain 
$\kappa_0\ge\hdots\ge \kappa_{2n}\ge 0$.

\section{Proof of Lemma \ref{lm:indicetot}}\label{s:IndiceTot}
\begin{lemma}
Assume that the determinant $\mathrm{det}\,B$ is of class $\mathcal{C}^1$ on $\bd$. Then the Maslov index $\kappa$ of $B$ is given by
$$\kappa=\mathrm{Ind}_{\mathrm{det}B(\bd)}(0)=\frac{1}{2\pi i}\int_{\bd}\frac{(\mathrm{det}B)'(\zeta)}{\mathrm{det}B(\zeta)}\,\mathrm{d}\zeta.$$
\end{lemma}

\begin{proof}
Assume that the partial indices come from the following factorization:
$$\forall\theta,\ B(e^{i\theta})=B^+(e^{i\theta})\left(\begin{array}{ccc}e^{i\kappa_0\theta}& &(0) \\ &\ddots& \\ (0)& & e^{i\kappa_{2n}\theta}\end{array}\right)B^{-}(e^{i\theta})$$
where $B^+$ extends holomorphically to $\D$ in some invertible matrix, and $B^-$ extends anti-holomorphically to $\D$ in some invertible matrix (that is, there exists $\tilde{B}^-$ holomorphic on $\D$ such that $B^{-}(\zeta)=\tilde{B}^-(1/\zeta)$ for all $\zeta\in\hat{\C}\setminus\D$).

Pick $0<r<1$ and set
$b_r^+(\theta)=\mathrm{det}(B^+(re^{i\theta}))$, $b_r^-(\theta)=\mathrm{det}(\tilde{B}^-(\overline{re^{i\theta}}))=\mathrm{det}(\tilde{B}^-(re^{-i\theta}))$ and $\beta_r(\theta)=b_r^+(\theta)r^\kappa e^{i\kappa\theta}b_r^-(\theta)$.
The path $\gamma_r=\beta_r([0;2\pi])$ does not meet 0, which allows us to define the index 
$$2\pi i\,\mathrm{Ind}_{\gamma_r}(0)=\int_{\gamma_r}\frac{\mathrm{d}\zeta}{\zeta}=\int_0^{2\pi}\frac{{b_r^+}'(\theta)}{b_r^+(\theta)}\mathrm{d}\theta+\int_0^{2\pi}i\kappa\mathrm{d}\theta+\int_0^{2\pi}\frac{{b_r^-}'(\theta)}{b_r^-(\theta)}\mathrm{d}\theta.$$
The index $\int_0^{2\pi}\frac{{b_r^+}'(\theta)}{b_r^+(\theta)}\mathrm{d}\theta$ is equal to the number of zeroes minus the number of poles of the holomorphic function $\mathrm{det}(B^+)$ in $r\bar{\D}$, that is, 0.
We also get $\int_0^{2\pi}\frac{{b_r^-}'(\theta)}{b_r^-(\theta)}\mathrm{d}\theta=0$, thus
$\mathrm{Ind}_{\gamma_r}(0)=\kappa$ for all $0<r<1$.
The two closed paths $\gamma_{1/2}$ and $\gamma_1$ have the same index since they are homotopic, thus:
$$\kappa=\mathrm{Ind}_{\gamma_{1/2}}(0)=\mathrm{Ind}_{\gamma_{1}}(0)=\frac{1}{2\pi i}\int_{0}^{2\pi}\frac{\beta'_1(\theta)}{\beta_1(\theta)}\mathrm{d}\theta,$$
which concludes the proof.
\end{proof}

Let us apply this statement to the matrix $B(\zeta)$ introduced in (\ref{B}):
$$\mathrm{det}\,B(\zeta)=1\times\begin{array}{|cc|}0&\zeta^2\\ \zeta^2&0\end{array}\times{\begin{array}{|cc|}O&\zeta\\ \zeta&0\end{array}\,}^{n-1}=(-1)^n\zeta^{2n+2},$$
which defines a $\mathcal{C}^1$-map on $\bd$. As a corollary, we obtain Lemma \ref{lm:indicetot}.

\bibliographystyle{amsplain}

\end{document}